\documentclass[a4paper,11pt]{extarticle}

\usepackage[T1]{fontenc}	
\usepackage[latin1]{inputenc}	

\usepackage{geometry}		
\usepackage{fancyhdr}		
\usepackage{setspace}		
\usepackage{lscape}		
\usepackage{multicol}		
\usepackage{makeidx}		
\usepackage{varioref}	
\usepackage{float, caption}

\usepackage{calc}		
\usepackage{ifthen}		
\usepackage{xspace}		

\usepackage{pifont}		
\usepackage{eurosym}		
\usepackage{soul}		
\usepackage{enumerate}		
\usepackage{verbatim}		
\usepackage{moreverb}
\usepackage{fancybox}
\usepackage{boxedminipage}		

\usepackage{array}		
\usepackage{multirow}		
\usepackage{tabularx}		
\usepackage{longtable}		

\usepackage{graphicx}		
\usepackage{epic}		
\usepackage{eepic}		
\usepackage{afterpage}		
\usepackage{rotating}		
\usepackage{pgf,tikz}
\usepackage{mathrsfs}
\usetikzlibrary{arrows}

\usepackage{amsmath}		
\usepackage{amssymb}		
\usepackage{mathrsfs}		
\usepackage[thmmarks,amsmath]{ntheorem}		
\usepackage{dsfont}

\newtheorem{prop}{Proposition}[subsection]

\newtheorem{theo}[prop]{Theorem}

\newtheorem{lem}[prop]{Lemma}
\newtheorem{ex}[prop]{Example}
\newtheorem{conj}[prop]{Conjecture}

\newtheorem{rem}[prop]{Remark}
\newtheorem*{theononnum}{Theorem}

\title{A family of centered random walks on weight lattices conditioned to stay in Weyl chambers}
\author{Vivien Despax\thanks{Laboratoire de Math\'ematiques et Physique Th\'eorique, UMR-CNRS 7350,
F\'ed\'eration Denis Poisson, FR-CNRS 2964,
Universit\'e Fran\c{c}ois--Rabelais de Tours, Parc de Grandmont, 37200 Tours, France, vivien.despax@lmpt.univ-tours.fr}}

\begin{document}
\maketitle
\begin{abstract} 
Under a natural asumption on the drift, the law of the simple random walk on the multidimensional first quadrant conditioned to always stay in the first octant was obtained by O'Connell in~\cite{OC2}. It coincides with that of the image of the simple random walk under the multidimensional Pitman transform and can be expressed in terms of specializations of Schur functions.
This result has been generalized in~\cite{LLP} and~\cite{LLP3} for a large class of random walks on weight lattices defined from representations of Kac-Moody algebras and their conditionings to always stay in Weyl chambers. In these various works, the drift of the considered random walk is always assumed in the interior of the cone.
In this paper, we introduce for some zero drift random walks defined from minuscule representations a relevant notion of conditioning to stay in Weyl chambers and we compute their laws. Namely, we consider the conditioning for these walks to stay in these cones until an instant we let tend to infinity. We also prove that the laws so obtained can be recovered by letting the drift tend to zero in the transitions matrices obtained in~\cite{LLP}. We also conjecture our results remain true in the more general case of a drift in the frontier of the Weyl chamber.
\end{abstract}

\tableofcontents

\section{Introduction}
Since the seminal work of O'Connell, various links between conditioned random walks and representation theory have been explored. In~\cite{LLP}, the authors generalize O'Connell's result to a wider class of random walks defined from minuscule representations of simple Lie algebras over $\mathbb{C}$. According to the standard root systems classification, the main result of ~\cite{OC2} appears in~\cite{LLP} as related to the vector representation in type $A$. In general, the first octant should be replaced by the closure of a Weyl chamber associated to the considered root system. Also one has to replace the path transformation used by O'Connell --- which is based on the Robinson-Schensted-Knuth correspondence --- by another map (the generalized Pitman transform) introduced by Biane, Bougerol and O'Connell in~\cite{BBO}. It is defined in terms of Littelmann's path model --- or, equivalently, in terms of Kashiwara's crystal basis theory --- applied to tensor powers of the minuscule representation used to define the steps of the random walk considered. Concerning these topics in representation theory, we refer the reader to the articles~\cite{Lit1},~\cite{Lit2},~\cite{Lit3} and~\cite{Kashi},~\cite{KashiNaka}. Beside these algebraic aspects, probabilistic tools from the Martin boundary theory are also used in~\cite{OC2} and in~\cite{LLP}. Mention also the possibility to use the Gessel-Zeilberger reflection principle (a multi-dimensional generalization of the classical geometric argument introduced by Andr\'{e}\ to solve the two candidates ballot problem) to derive a simpler combinatorial proof in the minuscule case, see \cite{Chha}.

Continuing in this vein but in the more general context of the representation theory of Kac-Moody algebras, the same authors proposed in \cite{LLP3} a different approach which does not require this time any result from the Martin boundary theory and can also be used to obtain the conditioning of random paths defined from non-minuscule representations (one should then replace the random walks by continuous stochastic processes). 

The assumption that the drift of the considered random walk belongs to the (open) Weyl chambers is an essential feature of all these works. This indeed yields a positive probability that it always remains in the closed cone, which
makes its conditioning easy to define. Nevertheless, one can observe that the
expressions obtained for the transition matrice of the conditioned process still make sense without this assumption. This follows from the fact that the
generalized Pitman transform applied to our random walk always provides us
with Markov chain whatever its drift.

It is then natural to ask if the Markov chain obtained as the Pitman
transform of a random walk with drift outside the (open) Weyl chamber still
has a probabilistic interpretation in terms of conditioning. The purpose of
this article is to answer this question in the zero drift case. Our approach is the following : one can first consider the law of the centered random walk conditioned by the event "Until instant $n$, the random walk stays in the closed cone" for any positive integer $n$ and let then the instant $n$ tends to infinity. One can also start from the known laws corresponding to a drift in the (open) Weyl chamber and let the drift tends to zero. The main result of our paper consists in proving that these two notions of conditioning coincide. Our proofs use algebraic properties of the random walks we study together with a recent probabilistic theorem of Denisov and Wachtel on the exit time from cones of centered random walks, see~\cite{DW}. We also refer the interested reader to ~\cite{Desp} for an elementary study of the simple random walks on the integers only based on the reflection principle and the strong law of large numbers. This approach could give a flavor of our general problem but it seems unfortunately difficult to adapt these arguments in higher dimension.

The paper is organized as follows. Section~\ref{recall} is a
recollection on known results on representation theory and the related random walks we study. We choose to restrict ourselves to the case of random walks defined from minuscule representations because a more general treatment would impose to consider random continuous trajectories (random
Littelmann paths in fact) and thus would require a more general version of the
deep result of Denisov and Wachtel. We refer to~\cite{Bour},~\cite{Hum},~\cite{Ser},~\cite{Ful} for basics on root systems, weight lattices and representation theory. We also review some specific facts about tensor multiplicities and the random walks we consider that are important for our purpose. According to~\cite{LLP} and~\cite{LLP3}, we define particular probability distributions on trajectories in order to permit an algebraic treatment of the conditioning to always stay in the (closed) Weyl chamber. We develop a three-dimensional working example to illustrate our main definitions. The law of these conditioned random walks in the case of a drift in the open cone is then recalled in paragraph~\ref{llp} as stated in~\cite{LLP}. In Section~\ref{main}, we turn to the specific case of the zero drift where the probability that our random walk always stays in the (closed) Weyl chamber becomes zero. We introduce the conditioning to stay in the (closed) Weyl chamber until an instant $n$ and then study what happens when $n$ tends to infinity. We then prove our main theorem by using Denisov and Wachtel estimate on the tail distribution of the exit time from this cone. Let us give a first statement of our main result. The reader is referred to Theorem~\ref{despax} for a definitive version.  

\begin{theononnum}
The transition probabilities corresponding to the conditioning to stay in the (closed) Weyl chamber until an instant $n$ of a centered random walk with steps the weights of a minuscule representation converge when $n$ tends to infinity. The limits so obtained can be recovered by letting the drift tends to zero in the transition probabilities corresponding to the same conditioning for random walks with the same steps but with drifts in the open Weyl chamber :

\begin{center}
\scriptsize{
\begin{tabular}{cccc}
\begin{boxedminipage}[c]{6.5cm}Law of the random walk\newline with drift in the (open) Weyl chamber\newline conditioned to stay in the (closed) Weyl chamber\newline until instant $n$\end{boxedminipage}& $\xrightarrow[n\to\infty]{}$ & \begin{boxedminipage}[c]{7.5cm}Law of the random walk\newline with drift in the (open) Weyl chamber\newline conditioned to always stay in the (closed) Weyl chamber\end{boxedminipage}\\
&&\\

$\downarrow \textbf{d}=0$& &  $\downarrow \textbf{d}=0$\\

&&\\

\begin{boxedminipage}[c]{6.5cm}Law of the random walk\newline with zero drift\newline conditioned to stay in the (closed) Weyl chamber\newline until instant $n$ \end{boxedminipage}& $\xrightarrow[n\to\infty]{}$ &\begin{boxedminipage}[c]{7.5cm} Law of the random walk\newline with zero drift\newline conditioned to always stay in the (closed) Weyl chamber
\end{boxedminipage}\end{tabular}}

\end{center}

\end{theononnum}

\section{Background}
\label{recall}
\subsection{Weight lattices associated to root systems}
\label{weight}
Let $d$ be a positive integer and $\mathfrak{g}$ a simple Lie algebra over $\mathbb{C}$ of rank $d$. Denote by $\mathfrak{h}$ a Cartan subalgebra of $\mathfrak{g}$. The relative root system $R$ is realized as a finite subset of a $d$-dimensional subspace $E$ of a standard Euclidean space $\mathbb{R}^D$ whose standard basis and standard inner product are respectively denoted by $\mathcal{B}=\left(\varepsilon_1,\ldots,\varepsilon_D\right)$ and $\left\langle \,.\, , \, .\,\right\rangle$. The classification of root systems distinguishes  four infinite families $A_d,B_d,C_d,D_d$ (corresponding respectively to the classical Lie algebras $\mathfrak{sl}_{d+1}$ ($d\geq1$), $\mathfrak{so}_{2d+1}$ ($d\geq 3$), $\mathfrak{sp}_{2d}$ ($d\geq2$), $\mathfrak{so}_{2d}$ ($d\geq3$)) and five exceptional types $E_6,E_7,E_8,F_4,G_2$.

We write $W=W\left( R\right) $ for the Weyl group of the root system. Let us choose a base $S=\left( \alpha_1,\ldots,\alpha_d\right)$ of $R$ and denote by $R^+$ the related finite set of positive roots. The family of fundamental dominant weights $\Omega=\left( \omega_1,\ldots,\omega_d\right)$ corresponding to this choice of $S$ is then defined by
\[2\frac{\left\langle \omega_i,\alpha_j\right\rangle }{\left\langle \alpha_j,\alpha_j\right\rangle }=\delta_{i,j}\qquad 1\leq i,j\leq d.\]
The Weyl chamber $C=C\left( R\right) $ is then the open convex cone of $E$ such that \[C=\bigoplus_{i=1}^{d} \mathbb{R}_{>0} \omega_i.\] We write $\overline{C}=\oplus_{i=1}^{d} \mathbb{R}_{\geq 0} \omega_i$ and $\partial C=\overline{C}\setminus C$ for, respectively, its closure and its frontier in $E$. The weight lattice $P=P\left( R\right) $ is the subgroup of $E$ generated by $\left\lbrace\omega_1,\ldots,\omega_d\right\rbrace $ : \[P=\bigoplus_{i=1}^d \mathbb{Z} \omega_i.\]It can be identified with an integral sublattice of $\mathfrak{h}_{\mathbb{R}}^*$. The subset of the dominant weights $P^+=P^+\left( R\right) $ is \[P^+=P\cap \overline{C}=\bigoplus_{i=1}^d \mathbb{Z}_{\geq 0} \omega_i.\] Although the Lie algebra $\mathfrak{gl}_{d+1} \left( \mathbb{C}\right)$ is not simple, we will consider in type $A_d$ the weight lattice of $\mathfrak{gl}_{d+1} \left( \mathbb{C}\right)$ rather than that of the simple Lie algebra $\mathfrak{sl}_{d+1} \left( \mathbb{C}\right)$.
 \begin{ex}
\label{weB3weights}
For the Lie algebra $\mathfrak{g}=\mathfrak{so_7 \left(\mathbb{C} \right) }$ and its Cartan subalgebra of diagonal matrices $\mathfrak{h}$ (the $B_3$ case), we have :
\begin{itemize}
\item[\textbullet]Rank : $d=3$
\item[\textbullet]$D=3$
\item[\textbullet]Root system : $R=\left\lbrace \pm\varepsilon_1,\pm \varepsilon_2,\pm \varepsilon_3,\pm \varepsilon_1\pm\varepsilon_2,\pm\varepsilon_1\pm\varepsilon_3,\pm\varepsilon_2\pm\varepsilon_3 \right\rbrace$
\item[\textbullet]Weyl group : $W$ is the semidirect product of $\mathfrak{S}_3$, the symmetric group operating by permutation on the set $\left\lbrace \varepsilon_1,\varepsilon_2,\varepsilon_3\right\rbrace$, and $\left(\mathbb{Z}/2\mathbb{Z} \right)^3 $ acting on the same set by $\varepsilon_i \mapsto \left(\pm 1\right) _i \varepsilon_i$. The group $W$ is of order $2^3 3!=48$.  
\item[\textbullet]A base : $S=\left(\alpha_1,\alpha_2,\alpha_3 \right)= \left( \varepsilon_1-\varepsilon_2,\varepsilon_2-\varepsilon_3,\varepsilon_3 \right)$
\item[\textbullet]Corresponding positive roots : $R^+=\left\lbrace  \varepsilon_1,\varepsilon_2,\varepsilon_3, \varepsilon_1 \pm \varepsilon_2,\varepsilon_1\pm\varepsilon_3,\varepsilon_2\pm\varepsilon_3\right\rbrace$
\item[\textbullet]Corresponding fundamental weights : $\Omega=\left(\omega_1,\omega_2,\omega_3 \right)= \left( \varepsilon_1,\varepsilon_1+\varepsilon_2,\frac{1}{2}\left(\varepsilon_1+\varepsilon_2+\varepsilon_3 \right) \right)$
\item[\textbullet]Corresponding Weyl chamber : $C=\left\lbrace x_1 \varepsilon_1+x_2 \varepsilon_2+x_3 \varepsilon_3 : \left(x_1,x_2,x_3 \right)\in\mathbb{R}^3 \enskip x_1>x_2>x_3>0  \right\rbrace$
\item[\textbullet]Corresponding Weight lattice : $P=\left( \mathbb{Z}\varepsilon_1 \oplus \mathbb{Z}\varepsilon_2\oplus \mathbb{Z}\varepsilon_3\right)+\mathbb{Z}\left(\frac{1}{2}\left(\varepsilon_1+\varepsilon_2+\varepsilon_3 \right) \right)$
\item[\textbullet]Corresponding set of dominant weights :\[P^+=\left\lbrace \lambda_1\varepsilon_1+\lambda_2\varepsilon_2+\lambda_3\varepsilon_3:\left(\lambda_1,\lambda_2,\lambda_3 \right)\in \mathbb{Z}^3 +\mathbb{Z}\left(\frac{1}{2},\frac{1}{2},\frac{1}{2} \right)  \enskip \lambda_1\geq \lambda_2 \geq \lambda_3\geq 0 \right\rbrace\]
\end{itemize}
\end{ex}

\subsection{Finite-dimensional representations}
\label{rep}
Consider the settings of the previous paragraph.
\subsubsection{Weights, characters, decomposition into irreducibles and tensor multiplicities}
Recall that $P$ is identified with an integral sublattice of $\mathfrak{h}_{\mathbb{R}}^*$. Given a finite-dimensional representation $V$ of $\mathfrak{g}$ and $\gamma$ in $P$, we denote by $V_\gamma$ the subspace of $V$ defined by \[V_\gamma =\left\lbrace v\in V : hv=\gamma\left(h \right)v \textit{ for any } h \textit{ in }\mathfrak{h}\right\rbrace. \]
When $V_\gamma$ is nonzero, we say that $\gamma$ is a weight of $V$ and the set of weights of the representation $V$ is denoted by $P\left( V\right) $. One has the decomposition \[V=\bigoplus_{\gamma\in P} V_\gamma=\bigoplus_{\gamma\in P\left( V\right) } V_\gamma.\]

The character of the finite-dimensional representation $V$ is then defined as the Laurent polynomial in $D$ variables $x=\left(x_1,\ldots,x_D\right)$ \[\mathrm{char}\,V\left( x\right) =\sum_{\gamma \in P }\dim V_\gamma \,x^\gamma=\sum_{\gamma \in P\left( V\right)  }\dim V_\gamma \,x^\gamma\]
where $x^\gamma =x_{1}^{\gamma_1} \ldots x_{D} ^{\gamma_N}$ if $\gamma=\sum_{i=1}^D \gamma_i \varepsilon_i$. Remark that \[\dim V =\mathrm{char}\,V\left( 1^D\right)\]
where we denote by $1^D$ the vector $\sum_{i=1}^D \varepsilon_i$.

The irreducible finite-dimensional representations of $\mathfrak{g}$ are labelled by the elements of $P^+$, that is by the dominant weights : to each dominant weight $\lambda$ corresponds a unique (up to isomorphism) irreducible finite dimensional representation of $\mathfrak{g}$ denoted by $V\left(\lambda\right)$. One has $V\left( 0\right) =\left\lbrace 0\right\rbrace $. Given a dominant weight $\lambda$ and a weight $\gamma$, the dimension of the subspace \[V\left( \lambda\right)_\gamma =\left\lbrace v\in V\left( \lambda\right)  : hv=\gamma\left(h \right)v \textit{ for any } h\textit{ in }\mathfrak{h}\right\rbrace\] is denoted by $K_{\lambda,\gamma}$. The finite set of weights of $V\left( \lambda\right) $ is denoted by $P\left( \lambda\right)$ instead of $P\left(V\left(  \lambda\right) \right)$ and we write $s_\lambda$ for the character of $V\left(\lambda \right)$. We thus have \[s_\lambda \left( x\right) =\mathrm{char}\,V\left( \lambda\right)\left( x\right) =\sum_{\gamma\in P}K_{\lambda,\gamma}x^\gamma=\sum_{\gamma\in P\left( \lambda\right) }K_{\lambda,\gamma}x^\gamma.\]
The family $\left(s_\lambda\right) _{\lambda \in P^+}$ is a basis of the vector space $\mathbb{C}^{W}\left[P\right]$. Any finite-dimensional representation $V$ of $\mathfrak{g}$ decomposes as a direct sum of irreducible components $V\left( \lambda\right) $ with finite multiplicities $m_{V,\lambda}$ :\[V\simeq \bigoplus_{\lambda \in  P^+} V\left( \lambda\right)^{\oplus m_{V,\lambda}}. \]

We use some standard notations when $V$ is obtained by tensor products of irreducible finite-dimensional representations. For any dominant weights $\lambda,\delta,\Lambda$ and any nonnegative integer $n$, the multiplicity of the irreducible component $V\left(\Lambda \right)$ in "the" decomposition into irreducibles of the finite-dimensional representation $V\left(\lambda\right)\otimes V\left(\delta\right)^{\otimes n}$ will be denoted by $f_{\Lambda / \lambda,\delta}^{n}$ :
\[V\left(\lambda\right)\otimes V\left(\delta\right)^{\otimes n}\simeq\bigoplus_{\Lambda\in P^+} V\left(\Lambda \right)^{\oplus f_{\Lambda/\lambda,\delta}^n}.\]
When $n=1$, we write $m_{\lambda,\delta}^{\Lambda}$ instead of $f_{\Lambda / \lambda,\delta}^{1}$. Remark that $m_{\lambda,\delta}^{\Lambda}=m_{\delta,\lambda}^{\Lambda}$.

In the sequel, the dominant weight $\delta$ will be fixed and we will then write $\lambda\leadsto\Lambda$ for two dominant weights $\lambda,\Lambda$ if $V\left( \Lambda\right)$ appears with a positive multiplicity in "the" decomposition of $V\left(\lambda \right)\otimes V\left(\delta \right)$ into irreducibles :
 \[V\left(\lambda \right)\otimes V\left(\delta \right)\simeq \bigoplus_{\Lambda\in P^+ }V\left(\Lambda \right)^{\oplus m_{\lambda,\delta}^{\Lambda}}=\bigoplus_{\lambda\leadsto\Lambda} V\left(\Lambda \right)^{\oplus m_{\lambda,\delta}^{\Lambda}}.\]
 The inequality \begin{equation}m_{\lambda,\delta}^{\Lambda}\leq K_{\delta, \Lambda -\lambda} \label{cryst}\end{equation}
 can be proved, for example, by using Kashiwara's crystal basis theory, see Proposition 5.3 in~\cite{LLP}.
 
\subsubsection{Minuscule weights and representations} 
If $\delta$ is a dominant weight such that its orbit under the action of $W$ contains $P\left( \delta\right)$, we say that $V\left( \delta\right)$ and $\delta$ are minuscule. Here is the type-classification of the minuscule weights.
\[\begin{tabular}{cccc}
\hline
Type & Minuscule weights & $D$ & Decomposition on $\mathcal{B}$ \\
\hline
$A_d$ & $\omega_i, i=1,\ldots,d$ & $d+1$ & $\varepsilon_1+\ldots+\varepsilon_i, i=1,\ldots,d$ \\

$B_d$ & $\omega_d$ & $d$ & $\frac{1}{2}\left(\varepsilon_1+\ldots+\varepsilon_d \right)$ \\

$C_d$ & $\omega_1$ & $d$ & $\varepsilon_1$ \\

$D_d$ & $\omega_1,\omega_{d-1},\omega_d$ & $d$ & $\varepsilon_1,\frac{1}{2}\left(\varepsilon_1+\varepsilon_2+\ldots-\varepsilon_d \right),\frac{1}{2}\left(\varepsilon_1+\varepsilon_2+\ldots+\varepsilon_d \right)$ \\

$E_6$ & $\omega_1,\omega_6$ & 8 & $\frac{2}{3}\left(\varepsilon_8-\varepsilon_7-\varepsilon_6 \right),\frac{1}{3}\left(\varepsilon_8-\varepsilon_7-\varepsilon_6 \right)+\varepsilon_5$\\

$E_7$ & $\omega_7$ & $8$ & $\varepsilon_6+\frac{1}{2}\left(\varepsilon_8-\varepsilon_7 \right)$   \\
\hline
\end{tabular}\]

Minuscule weights are fundamental weights and, in the minsucule case, every subspace\[V\left( \delta\right)_\gamma =\left\lbrace v\in V\left( \delta\right)  : hv=\gamma\left(h \right)v \textit{ for any } h\textit{ in }\mathfrak{h}\right\rbrace\qquad \gamma\in P\]which is nonzero is one-dimensional and so \[s_\delta \left( x\right) =\sum_{\gamma\in P\left( \delta\right) } x^{\gamma}.\]When $\delta$ is minuscule, one can also prove \[m_{\lambda,\delta}^{\Lambda} = K_{\delta,\Lambda-\lambda}\in\left\lbrace 0,1\right\rbrace \qquad \lambda,\Lambda\in P^+.\]
In this case, we have then\[\lambda\leadsto \Lambda  \iff  m_{\lambda,\delta}^{\Lambda} =1 \iff  \Lambda-\lambda \in P\left( \delta\right) \qquad \lambda,\Lambda\in P^+,\]
so \[V\left(\lambda \right)\otimes V\left(\delta \right)=\bigoplus_{\lambda\leadsto\Lambda}V\left(\Lambda \right)\qquad \lambda \in P^+\]
and this gives \begin{equation}s_\lambda \left( x\right) s_\delta \left( x\right) =\sum_{\lambda\leadsto\Lambda}s_\Lambda \left(x\right)\qquad \lambda\in P^+.\label{pieri}\end{equation}

\begin{ex}
\label{weB3minus}
Consider the settings of Example~\ref{weB3weights}.
\begin{itemize}
\item[\textbullet]Minuscule weight : $\omega_3=\frac{1}{2}\left( \varepsilon_1+\varepsilon_2+\varepsilon_3\right) $
\item[\textbullet]Weights of the minuscule weight : $P\left( \omega_3\right) =\left\lbrace \frac{1}{2}\left( \pm \varepsilon_1\pm \varepsilon_2\pm\varepsilon_3\right)  \right\rbrace$
\item[\textbullet]Character of the minuscule weight : $s_{\omega_3} \left( x_1,x_2,x_3\right)=\sum x_1^{\pm \frac{1}{2}} x_2^{\pm \frac{1}{2}} x_3 ^{\pm \frac{1}{2}}$
\end{itemize}
\end{ex}

Given a minuscule weight $\delta$, tensor multiplicities have a combinatorial interpretation that will be important for our purpose. When $\gamma,\Gamma$ are in $P$ and when $n$ is a positive integer, we denote by $\pi_n \left( \gamma,\Gamma\right)$ the set \[\left\lbrace \left(\pi_0,\pi_1,\ldots,\pi_n \right)\in P^{n+1}:\left( \pi_0,\pi_n\right) =\left( \gamma,\Gamma\right) ,\pi_{i+1}-\pi_i\in P\left( \delta\right) \enskip 0\leq i \leq n-1  \right\rbrace,\] that is the set of paths in $P$ with steps in $P\left( \delta\right) $ of length $n$ starting at $\gamma$ and ending at $\Gamma$. When $\lambda,\Lambda$ are two dominant weights, the subset of paths always remaining in $P^+$ is denoted by $\pi_{n}^+ \left( \lambda,\Lambda\right)$. Again, the following result can be proved thanks to considerations on Kashiwara crystals or Littelmann paths, see Corollary 5.5 in~\cite{LLP}.

\begin{prop}
\label{combi}
We have\[\mathrm{Card}\, \pi_{n}^+\left(\lambda,\Lambda  \right)=f_{\Lambda/\lambda,\delta}^n\qquad \lambda,\Lambda\in P^+,\,n\geq 1.\]
\end{prop}

\subsection{A family of random walks on the weight lattice}
\label{family}
Consider the settings of paragraph~\ref{weight} and let us choose once for all a minuscule weight $\delta$. We endow the finite set $P\left( \delta\right)$ with probability distributions that have the property that the probability of any path in $P$ with steps in $P\left( \delta\right)$ only depends on its two extremal points and its length. They are constructed as follow. Fix $\theta=\left(\theta_1,\ldots,\theta_d \right)$ in $\mathbb{R}_{>0}^d$. The dominant weight $\delta$ being minuscule, given a step $s$ in $P\left( \delta\right)$, there exists $w$ in the Weyl group $W$ such that $s=w\delta$. But, $\delta$ being a dominant weight, $\delta-w\delta$ is a positive root : we can then write $\delta-s=\delta-w\delta=\sum_{i=1}^d m_i \alpha_i$ where the $m_i$ are nonnegative integers. Let us then denote by $\theta^{\left[ s\right] }$ the positive real $\theta^{\left[s\right] }=\theta_1^{m_1}\ldots \theta_d^{m_d}$ and define \begin{equation}\Sigma=\sum_{s\in P\left( \delta\right) } \theta^{\left[s\right] }.\label{sum}\end{equation} The probability $p_s=p_{s}\left(\theta\right) $ of the step $s$ is then defined by
\[p_{s}=\frac{\theta^{\left[s\right]} }{\Sigma }.\]
With this construction, one can remark that if $s,S$ in $P\left( \delta\right)$ are such that $S=s-\alpha_i$ (geometrically, $S$ is the image of $s$ under the orthogonal reflexion of $E$ of hyperplan $\left\lbrace \alpha\right\rbrace ^{\perp} $) for some $i$ in $\left\lbrace 1,\ldots,d\right\rbrace $, then we have \[p_{s}\times\theta_i=p_{S}.\] We can represent this situation by the following oriented and colored graph \[\begin{array}{ccc}\boxed{s}&\xrightarrow{i}&\boxed{S}\end{array}.\]

In order to obtain a formulation of our result in terms of characters, we introduce $x=x\left( \theta\right)=\sum_{i=1}^D x_i \varepsilon_i$, a solution in $\oplus_{i=1}^D\mathbb{R}_{>0} \varepsilon_i$ of the system
\begin{equation}x^{\alpha_i}=\theta_i ^{-1}\qquad 1\leq i\leq d.\label{syst}\end{equation}
By using the Log function and elementary linear algebra, observe that this system possesses at least one solution since the family $S=(\alpha_1,\ldots,\alpha_d)$ of $\mathbb{R}^D$ is of rank $d$.
We then can rewrite the sum $\Sigma$
\[\Sigma=x^{-\delta} s_{\delta} \left( x\right)\]
and the probability distribution $\left( p_{s}\right) _{s\in P\left( \delta\right)  }$ is now given by
\[p_s =\frac{x^s}{s_\delta \left(x\right) }\qquad s\in P\left( \delta\right) .\]
Also notice that when one chooses $\theta=\left( 1,\ldots,1\right)=1^d$, one has the obvious solution $x=\sum_{i=1}^D \varepsilon_i=1^D$ and this yields the uniform distribution on $P\left( \delta\right)$ : in this particular case, each step has probability ${1}/{\dim \,V\left(\delta \right)  }$.

For the remainder of this paragraph, we fix a $\theta$ in $\mathbb{R}_{>0}^d$ and $\left( p_{s}\right) _{s\in P\left( \delta\right)  }$, the associated probability distribution on $P\left( \delta\right) $ constructed this way. We also pick $x$ in $\oplus_{i=1}^D\mathbb{R}_{>0}\varepsilon_i$, a solution of System~\eqref{syst}. We endow the denumerable set $P$ with the $\sigma$-algebra of all subsets. Let $\nu_0$ be a probability measure on $P$ with full support\[\nu_0\left( \gamma\right) >0\qquad \gamma\in P.\] Consider $\left( \mathcal{S}\left( 0\right) ,\mathcal{X}\left( 1\right) ,\mathcal{X}\left( 2\right), \ldots \right) $ a sequence of independent random variables defined on a probability space $\left(\Omega,\mathcal{T},\mathbb{P}\right)$ such that 
\begin{enumerate}
\item $\mathcal{S}\left(0\right) : \Omega\rightarrow P$ has law $\nu_0$.
For any $\gamma$ in $P$, we write $\mathbb{P}_\gamma$ instead of $\mathbb{P}_{\mathcal{S}\left( 0\right) =\gamma}$ or $\mathbb{P}\left[ \,.\,|\,\mathcal{S}\left( 0\right) =\gamma\right] $ for the conditional probability defined by\[\mathbb{P}_\gamma \left[ \,.\,\right] =\frac{\mathbb{P}\left[ \,.\, \cap \mathcal{S}\left( 0\right) =\gamma\right] }{\mathbb{P}\left[\mathcal{S}\left( 0\right) =\gamma \right] }.\]
\item $\mathcal{X}=\left( \mathcal{X}\left( n\right):\Omega\rightarrow P \right)_{n\geq 1}$ are identically distributed with law $\nu=\sum_{s \in P\left( \delta\right) }p_s \delta_{\left\lbrace s\right\rbrace }$. We denote by \[\textbf{d}=\textbf{d}\left( \theta\right) =\mathbb{E}\left[  \mathcal{X}\left( 1\right) \right] =\sum_{s\in P\left( \delta\right) }\mathbb{P}\left[ \mathcal{X}\left( 1\right) =s\right] =\frac{1}{s_{\delta} \left( x\right) }\sum_{s\in P\left( \delta\right) }x^s s\] the common expectation of the random variables $\mathcal{X}$.  
\end{enumerate}
We define a sequence of random variables $\mathcal{S}=\left(\mathcal{S}\left( n\right)  \right)_{n\geq 1}$ by setting\[\mathcal{S}\left(n\right)=\sum_{i=1}^n \mathcal{X}\left( i\right) \qquad n\geq 1.\] The sequence of random variables $\mathcal{W}=\left(\mathcal{W}\left( n\right)  \right)_{n\geq 0} $ defined by $\mathcal{W}\left( 0\right) =\mathcal{S}\left( 0\right) $ and\[\mathcal{W}\left( n\right) =\mathcal{S}\left( 0\right) +\mathcal{S}\left( n\right) \qquad n\geq 1\] is then called the random walk on $P$ with steps in $P\left( \delta\right)$, with initial distribution $\nu_0$ and with law of increments $\nu$. Its drift is the vector $\textbf{d}$. The sequence $\mathcal{W}$ is a (time-homogeneous) Markov chain with initial distribution $\nu_0$, state space $P$ and transition matrix $p=p\left( \theta\right)$ given by \[\begin{aligned}p\left( \gamma,\Gamma\right)&=\mathbb{P}\left[\mathcal{W}\left(i+1\right)=\Gamma\,|\,\mathcal{W}\left(i\right) =\gamma \right]\\&=\nu\left( \Gamma-\gamma\right)\\&= \begin{cases} \frac{x^{\Gamma-\gamma}}{s_\delta \left( x\right) } & \text{if } \Gamma-\gamma\text{ is in }P\left( \delta\right)  \\ 0 & \text{otherwise}\end{cases}\end{aligned}\qquad \gamma,\Gamma\in P,\,i\geq 0.\]
With the previous construction, one can use independance to check that for any $\gamma,\Gamma$ in $P$ and any positive integer $n$, we have \begin{equation}\mathbb{P}_{\gamma}\left[ \mathcal{W}\left( n\right) =\Gamma\right]=\mathbb{P}\left[\bigcap_{i=1}^n \lambda+\mathcal{S}\left( i\right)=\pi_i\right]=\frac{x^{\Gamma-\gamma}}{s_\delta \left( x\right)^n }\qquad \left(\pi_0,\pi_1,\ldots,\pi_n \right) \in \pi_n \left(\gamma,\Gamma\right).\label{centr}\end{equation}
 As desired, this number only depends on the two extremal points and on the length of the considered path.
 
 In Example~\ref{probaB3}, one can see that the drift $\textbf{{d}}$ is in the open cone $C$ if and only if the parameter $\theta$ is in $\left] 0,1\right[^3 $ and this vector is zero if and only if $\theta=1^3$, that is if and only if $\left( p_s\right) _{s\in P\left( \omega_3\right) }$ is the uniform distribution on $P\left(\omega_3\right)$. This is a general fact, but it is not a trivial one, see Lemma 7.2 in~\cite{LLP} for a proof.

\begin{ex}
\label{probaB3}
Consider the settings of Example~\ref{weB3minus}.
\begin{itemize}
\item[\textbullet]Solution of System~\eqref{syst} : $\left( x_1,x_2,x_3\right) =\left(\frac{1}{\theta_1 \theta_2 \theta_3},\frac{1}{\theta_2 \theta_3},\frac{1}{\theta_3}\right)$
\item[\textbullet]Sum~\eqref{sum} : \[\Sigma=1+\theta_3+\theta_2 \theta_3 +\theta_1\theta_2\theta_3+\theta_2 \theta_3^2+\theta_1 \theta_2 \theta_3 ^2+\theta^1\theta_2^2 \theta_3^2+\theta_1 \theta_2^2\theta_3^3\]
\item[\textbullet]Relation between the steps :  \[\scriptsize{\begin{array}{ccccccccccccc}
&&&&&          &\boxed{-++} &&&&&&\\
&&&&& \overset{1}{\nearrow} &  & \overset{3}{\searrow} &&&&&\\
\boxed{+++} & \xrightarrow{3}  & \boxed{++-} & 
\xrightarrow{2} & \boxed{+-+} & & & &\boxed{-+-} & \xrightarrow{2}  &\boxed{--+} & \xrightarrow{3} &  \boxed{---}\\
&&&&& \underset{3}{\searrow} &  & \underset{1}{\nearrow} &&&&&\\ 
&&&&&          &\boxed{+--} &&&&&& \end{array}}\]
Here, the eight step-vectors $\frac{1}{2}\left(\pm \varepsilon_1\pm \varepsilon_2\pm \varepsilon_3 \right)$ are replaced by the signs of their coordinates in the standard basis $\mathcal{B}$.
\item[\textbullet]Probability distribution $p$ on $P\left( \omega_3\right)$ :
\begin{itemize}
\item[\textbullet]$p_{+++}= \frac{1}{\Sigma}$
\item[\textbullet]$p_{++-}= \theta_3 p_{+++}= \frac{\theta_3}{\Sigma}$
\item[\textbullet]$p_{+-+} = \theta_2 p_{++-}= \frac{\theta_2\theta_3}{\Sigma}$
\item[\textbullet]$p_{-++}= \theta_1 p_{+-+}=\frac{\theta_1 \theta_2 \theta_3}{\Sigma}$
\item[\textbullet]$p_{+--}= \theta_3 p_{+-+}=\frac{\theta_2  \theta_3 ^2}{\Sigma}$
\item[\textbullet]$p_{-+-}= \theta_1 p_{+--}=\theta_3 p_{-++}=\frac{\theta_1 \theta_2 \theta_3^2}{\Sigma}$
\item[\textbullet]$p_{--+}=\theta_2 p_{-+-}=\frac{\theta_1 \theta_2^2  \theta_3 ^2}{\Sigma}$
\item[\textbullet]$p_{---}=\theta_3 p_{--+}=\frac{\theta_1  \theta_2^2\theta_3 ^3}{\Sigma}$
\end{itemize}
\item[\textbullet]Relations between probabilities :
\begin{itemize}
\item[\textbullet]$p_{++-} p_{+-+}=p_{+--}p_{+++}$
\item[\textbullet]$p_{++-}p_{-++}=p_{-+-}p_{+++}$
\item[\textbullet]$p_{++-}p_{--+}=p_{---}p_{+++}$
\item[\textbullet]$p_{+-+}p_{-++}=p_{--+}p_{+++}$
\item[\textbullet]$p_{+-+}p_{-+-}=p_{---}p_{+++}$
\item[\textbullet]$p_{-++}p_{+--}=p_{---}p_{+++}$
\end{itemize}
\item[\textbullet]Drift :\[\begin{aligned}\textbf{d}&=\frac{\left(1-\theta_1\right)\theta_2\theta_3 \left(1+\theta_3 \right) }{\Sigma}\omega_1+\frac{\left( 1-\theta_2\right)\theta_3\left( 1+\theta_1 \theta_2 \theta_3\right)  }{\Sigma}\omega_2\\&+\frac{\left(1-\theta_3 \right)\left(1+\theta_2\theta_3+\theta_1\theta_2\theta_3+\theta_1 \theta_2^2 \theta_3^2 \right)  }{\Sigma}\omega_3\end{aligned}\]   
\end{itemize}
\end{ex}

Let $P'$ be a denumerable set and $p'$ a stochastic matrix on $P'$, that is a map $p':P'\times P'\rightarrow \left[ 0,1\right]$ such that $\sum_{y\in P'}p'\left( x,y\right) =1$ for any $x$ in $P'$. One can build a measurable space $\left(\Omega',\mathcal{T}'\right)$ such that, for each probability measure $\nu_0'$ on $\left(\Omega',\mathcal{T}'\right)$, there exists a probability measure $\mathbb{P}_{\nu_0'}$ on $\left(\Omega',\mathcal{T}'\right)$ and a sequence of $P'$-valued maps defined on this probability space which is a Markov chain with initial distribution $\nu_0'$, state space $P'$ and transition matrix $p'$. For a detailed construction, see~\cite{W} for example. We adopt this point of view to define a Markov chain from $\mathcal{W}$ that we call "the" random walk conditioned to always stay in the dominant weights.

\subsection{Conditioning to always stay in the (closed) Weyl chamber when the drift is in the (open) Weyl chamber}
\label{llp}
Let us consider the event\[\mathcal{W}\in\overline{C}=\bigcap_{n\geq 0}\mathcal{W}\left( n\right)\in\overline{C}=\bigcap_{n\geq 0}\mathcal{W}\left( n\right)\in P^+. \] For any $\lambda$ in $P^+$, we also denote by $\lambda+\mathcal{S}\in \overline{C}$ the event \[\lambda+\mathcal{S}\in\overline{C}=\bigcap_{n\geq 1} \lambda+\mathcal{S}\left( n\right)  \in \overline{C}=\bigcap_{n\geq 1} \lambda+\mathcal{S}\left( n\right)  \in P^+.\]The set $\overline{C}$ being stable under addition, we have the inclusions\[\mathcal{S}\in\overline{C}\subset \lambda+\mathcal{S}\in\overline{C} \qquad \lambda\in P^+.\] Using independence and the previous remark, we have thus \begin{equation}\mathbb{P}_{\lambda}\left[\mathcal{W}\in\overline{C} \right]=\mathbb{P}\left[\lambda+\mathcal{S}\in\overline{C} \right]\geq \mathbb{P}\left[\mathcal{S}\in\overline{C}\right]=\mathbb{P}_0 \left[ \mathcal{W}\in\overline{C}\right]  \qquad \lambda\in P^+ \label{inc}\end{equation}
and the total probability formula gives then\begin{equation}\mathbb{P}\left[\mathcal{W}\in\overline{C} \right]=\sum_{\lambda\in P^+}\nu_0 \left(\lambda\right)\mathbb{P}\left[\lambda+\mathcal{S}\in\overline{C} \right].\label{tpf}\end{equation}

If the event $\mathcal{S}\in\overline{C}$ has a positive probability, so has the event $\mathcal{W}\in\overline{C}$ by~\eqref{inc} and~\eqref{tpf}. Then, from a probabilistic point of view, to take the positive probability event $\mathcal{W}\in \overline{C}$ into account leads to introduce the conditional probability $\mathbb{P}_{ \mathcal{W} \in \overline{C}}$ defined by \[\mathbb{P}_{\mathcal{W} \in \overline{C}}\left[ \,.\,\right]  =\frac{\mathbb{P}\left[ \,.\cap \mathcal{W} \in \overline{C} \right] }{\mathbb{P}\left[ \mathcal{W} \in \overline{C}\right] }.\]

Let us observe how the transition probability between two states chosen in $P^+$ is modified if one replaces $\mathbb{P}$ by $\mathbb{P}_{\mathcal{W} \in \overline{C} }$. One can easily derive the following proposition from~\eqref{finiteinfinite} with Proposition~\ref{condifinite} --- its finite-time equivalent --- or prove it directly.

\begin{prop}
\label{condpos}
Assume the event $\mathcal{S}\in\overline{C}$ has a positive probability. We have
\[\mathbb{P}_{\mathcal{W}\in \overline{C}} \left[ \mathcal{W}\left( i+1\right)=\Lambda \, |\, \mathcal{W}\left(i\right) =\lambda\right] =p\left( \lambda,\Lambda\right) \frac{\mathbb{P}_{\Lambda}\left[\mathcal{W}\in \overline{C} \right] }{\mathbb{P}_{\lambda}\left[\mathcal{W} \in \overline{C} \right] }\qquad \lambda,\Lambda\in P^+,\,i\geq 0.\]
\end{prop}
This proposition leads us to introduce the function\[\begin{array}{ccccc}
   h & : & P^+ & \longrightarrow &  \mathbb{R}\\
    & & \lambda & \longmapsto & \mathbb{P}_{\lambda}\left[\mathcal{W}\in\overline{C} \right]= \mathbb{P}\left[ \lambda+\mathcal{S}\in\overline{C}\right] 
   \label{hfunc}\end{array}.\]
Observe that the transitions of our random walk under $\mathbb{P}_{\mathcal{W}\in \overline{C}}$ are obtained by a particular transformation of the restriction of the matrix $p$ to the subset $P^+$. In the general setting, given a pair $\left(P',p' \right) $ where $P'$ is a denumerable set and $p'$ is a substochastic matrix on $P'$, that is a map $p':P'\times P'\rightarrow \left[ 0,1\right]  $ such that $\sum_{y\in P'}p'\left( x,y\right) \leq 1$ for any $x$ in $P'$, and a positive function $h':P'\to \mathbb{R}_{>0}$, we can define another matrix $p'_{h'}$ on $P'$ by setting\[p'_{h'} \left( x,y\right)=p'\left( x,y\right)\frac{h'\left( y\right) }{h'\left( x\right) }\qquad x,y\in P'.\] Then the matrix $p'_{h'}$ is stochastic on $P'$ if and only if $h'$ satisfies \[h'\left( x\right) =\sum_{y\in P'}p'\left( x,y\right) h'\left( y\right) \qquad x\in P'.\]In this case, we say that $h'$ is a positive harmonic function on $P'$ for $p'$. As we will see it in the following, in the positive drift case, the function $h$ that we introduced previously is an example of positive harmonic function on $P^+$ for the restriction of the matrix $p$ to $P^+$. Then, Proposition~\ref{condpos} shows that the map $p^{+}=p^{+}\left( \theta\right)$ defined on $P^+\times P^+$ by\[p^{+}\left( \lambda,\Lambda\right)=p\left( \lambda,\Lambda\right) \frac{h\left(\Lambda\right)}{h\left( \lambda\right) }\qquad \lambda,\Lambda\in P^+\label{condiposdrift}\]
is a stochastic matrix on $P^+$. Any Markov chain associated to the pair $\left( P^+,p^{+}\right)$ is called "the" random walk with law $\nu$ conditioned to always stay in the dominants weights and we denote it by $\mathcal{W}^+$.

\begin{prop}
For any drift $\textbf{d}$, the function $h$ is harmonic on the set $P^+$ for the restriction of the matrix $p$ to $P^+$.
\end{prop}
\textit{Proof.} The proof is almost direct. It is not more difficult than in the most elementary case of simple random walks on the integers, see~\cite{Desp}. $\Box$\\

Under the crucial assumption that the drift of the random walk $\mathcal{W}$ is in the open convex cone $C$, Lecouvey, Lesigne and Peigné achieved to explicit the law of $\mathcal{W}^+$ in terms of irreducibles characters, see Theorem 7.6 and Corollary 7.7 in~\cite{LLP}.

\begin{theo}[Lecouvey, Lesigne, Peigné]
\label{llptheo}	
Assume that the drift $\mathrm{\textbf{d}}$ is in the (open) Weyl chamber $C$ (equivalently : the parameter $\theta$ is in $\left]0,1\right[^d$). We have\[h\left( \lambda\right) =x^{-\lambda}s_\lambda \left( x\right) \prod_{\alpha\in R^+}\left( 1-x^{-\alpha}\right) \qquad \lambda\in P^+.\]The transition matrix $p^+$ of the Markov chain $\mathcal{W}^+$ is given by
\[p^+\left(\lambda,\Lambda \right)=p \left( \lambda,\Lambda\right) \frac{x^{-\Lambda} s_\Lambda \left( x\right) }{x^{-\lambda} s_\lambda \left( x\right)}=\begin{cases}\frac{1}{s_\delta \left( x \right) }\frac{s_{\Lambda}\left( x \right) }{s_{\lambda}\left(x \right)}&\text{if }\lambda\leadsto \Lambda \\ 0 & \text{otherwise}\end{cases}\qquad \lambda,\Lambda\in P^+.\]
\end{theo}
As we noticed it in the introduction, these expressions still make sense when the parameter $\theta$ is $1^d$. In the next section, we prove that the equality obtained for this value of the parameter can be still probabilistically interpreted in terms of conditioning for the corresponding centered random walk to always stay in the (closed) Weyl chamber. 

\section{The zero drift case}
\label{main}
\subsection{Conditioning to stay in the dominant weights until an instant}
We introduce $\left(\mathcal{W}\left( \leq n\right)  \right)_{n\geq 0}$ the increasing sequence of positive probability events defined by\[\mathcal{W}\left( {\leq n}\right) \in \overline{C}=\bigcap_{i=0}^n \mathcal{W}\left( {i}\right)  \in \overline{C}=\bigcap_{i=0}^n \mathcal{W}\left( {i}\right)  \in P^+ \qquad n\geq 0\] which limit is $\mathcal{W}\in\overline{C}$. We also introduce for each $\lambda$ in $P^+$ the increasing sequence $\left(\lambda+\mathcal{S}\left( \leq n\right)  \right)_{n\geq 1}$ defined by \[\lambda+\mathcal{S}\left( {\leq n}\right) \in \overline{C}=\bigcap_{i=1}^n  \lambda+\mathcal{S}\left( {i}\right)  \in \overline{C}=\bigcap_{i=1}^n \lambda+\mathcal{S}\left( {i}\right)  \in P^{+}\qquad n\geq 1\]
which limit is $\lambda+\mathcal{S}\in\overline{C}$.  Observe that we have \[\mathbb{P}\left[ \mathcal{W}\left( \leq n\right) \right] =\sum_{\lambda\in P^+}\nu_0 \left( \lambda\right) \mathbb{P}_{\lambda}\left[ \mathcal{W}\left( \leq n\right) \right]\qquad n\geq 0,\] \[\mathbb{P}_{\lambda}\left[\mathcal{W}\left(\leq n \right)  \right]=\mathbb{P}\left[ \lambda+\mathcal{S}\left( \leq n\right) \right]\qquad \lambda\in P^{+},\, n\geq 1,\]
\[\mathbb{P}\left[\lambda+\mathcal{S}\left( {\leq n}\right)  \in \overline{C}\right]\geq\mathbb{P}\left[ \mathcal{S}\left( {\leq n}\right)  \in\overline{C}\right]\geq\mathbb{P}\left[\bigcap_{i=1}^n \mathcal{X}\left( i\right) =+1  \right]=p_{\delta} ^n>0 \qquad  \lambda\in P^{+},\, n\geq 1.\]
When the drift is zero, Theorem~\ref{dw} shows in particular that the event $\mathcal{W}\in\overline{C}$ has probability zero. Then we can use these positive probability approximations of the event $\mathcal{W}\in\overline{C}$ in order to build a notion of conditioning of always staying in the dominant weights in the zero drift case. Let us introduce the sequence of conditional probabilities $\left( \mathbb{P}_{\mathcal{W}\left( {\leq n}\right)  \in {\overline{C}}} \right)_{n\geq 0}$ defined by \[\mathbb{P}_{\mathcal{W}\left( {\leq n}\right)  \in {\overline{C}}} \left[ \,.\,\right]=\frac{\mathbb{P}\left[\,.\cap \mathcal{W}\left( {\leq n}\right) \in \overline{C} \right]}{\mathbb{P}\left[\mathcal{W}\left( {\leq n}\right) \in \overline{C}\right]} \qquad n\geq 0.\]

Since the sequence $\left(\mathbb{P}\left[ \mathcal{W}\left( {\leq n}\right) \in\overline{C}\right]  \right) _{n\geq 0} $ converges to $\mathbb{P}\left[ \mathcal{W}\in\overline{C}\right] $, one has, under the assumption that $\mathbb{P}\left[\mathcal{W}\in\overline{C}\right]$ is positive (that is when the drift is in the (open) Weyl chamber $C$ or when $\theta$ is in $\left] 0,1\right[^d$),\[\mathbb{P}_{ \mathcal{W}\left( {\leq n}\right)  \in {\overline{C}}} \left[E \right]\xrightarrow[n \to \infty]{}\mathbb{P}_{\mathcal{W}\in\overline{C}}\left[ E\right] \qquad E\in\mathcal{T}.\] If we set \[p_{n}^{+}\left(\lambda,\Lambda,i\right)=\mathbb{P}_{\mathcal{W}\left( {\leq n}\right)  \in {\overline{C}}} \left[\mathcal{W}\left( {i+1}\right)=\Lambda \,|\,\mathcal{W}\left( i\right)=\lambda\right]\qquad \lambda,\Lambda\in P^+,\,n\geq 0,\, i\geq 0,\]the previous remark gives in particular
\begin{equation}p_{n}^{+}\left( \lambda,\Lambda,i\right) \xrightarrow[n\to\infty]{}p^+ \left( \lambda,\Lambda\right) \qquad \lambda,\Lambda\in P^+,\,i\geq 0 \label{finiteinfinite}\end{equation}
whenever the parameter $\theta$ is in $\left] 0,1\right[^d$.

When the drift is zero, that is when the simple random walk considered is the centered one, it is no longer obvious whether the sequences \[\left(\mathbb{P}_{\mathcal{W}\left( {\leq n}\right) \in\overline{C} } \left[ \mathcal{W}\left( {i+1}\right)=\Lambda \, |\,\mathcal{W}\left( {i}\right) =\lambda\right]\right)_{n\geq 0}\qquad \lambda,\Lambda\in P^+,\,i\geq 0\] are convergent, since the denominator $\mathbb{P}\left[ \mathcal{W}\left( {\leq n}\right) \in\overline{C}\right]$ tends to zero when $n$ tends to infinity.

As done in the positive drift case, we start by observing how the transition probability between two states in $P^+$ is modified when we replace $\mathbb{P}$ by one of the $\mathbb{P}_{\mathcal{W}\left( {\leq n}\right) \in\overline{C}}$.
\begin{prop}
\label{condifinite}
Let $\lambda,\Lambda$ be in $P^+$. One has\[p_{n}^{+}\left(\lambda,\Lambda,i\right) =\begin{cases} p\left( \lambda,\Lambda\right) \frac{\mathbb{P}_{\Lambda}\left[ \mathcal{W}\left( \leq{n-i-1}\right)\right]  } {\mathbb{P}_\lambda \left[ \mathcal{W}\left( \leq{n-i}\right)\right]}&\text{if } i\leq n-2\\p\left( \lambda,\Lambda\right) &\text{if } i> n-2 \end{cases} \qquad n\geq 2,\,i\geq 0.\]
In particular, one has for any nonnegative integer $i$\[p_{n}^{+}\left(\lambda,\Lambda,i\right)=p_{n-i}^{+}\left(\lambda,\Lambda,0\right) \qquad  n\geq i+2.\]
\end{prop}
\textit{Proof.} Once again, the proof is not more difficult than in the most elementary case of simple random walks on the integers, see~\cite{Desp}. $\Box$\\

\subsection{Exit time from a cone in the zero drift case}
  
 Given a dominant weight $\lambda$, we introduce the random variable \[\tau_\lambda =\inf \left\{n\geq 1:\lambda+\mathcal{S}\left( n\right)  \not\in \overline{C} \right\},\] that is the first time the random walk starting at $\mu$ exits from the closed cone $\overline{C}$. In a more general setting, Denisov and Wachtel found out an asymptotic for the tail distribution of the exit time from a cone of a centered random walk under some conditions on the cone and the moments of the law of increments, see~\cite{DW}. As $\overline{C}$ is a closed convex polyhedral cone (i.e. the convex hull of a finite set of half lines) and as our set of steps $P\left(\delta \right) $ is finite, the conditions they require are fulfilled in our particular situation. 
 
\begin{theo}[Denisov,Wachtel]
\label{dw}
Assume that the drift of the random walk $\mathcal{W}$ is zero, that is $\left( p_s\right)_{s\in P\left( \delta\right) }$ is the uniform distribution on the set of steps $P\left( \delta\right) $ (equivalently : $\theta$ is $1^d$). There exists two positive reals $\kappa_1,\kappa_2$ and a positive function $V:\overline {C} \rightarrow \mathbb{R}_{>0}$ such that \[\mathbb{P}
 \left(\tau_\lambda >n \right) \underset{n\to \infty}{\sim} \kappa_1 V\left(\lambda \right)n^{-\kappa_2} \qquad \lambda\in P^+.\]
\end{theo}

\begin{rem}
\begin{enumerate}
\item From this theorem, the existence of the positive function $V$ and the polynomial rate of convergence to zero are the only facts that we need in the following, see Proposition~\ref{harm}. For a complete presentation including detailed discussion about hypotheses, constants and the positive function appearing in the conclusion of this theorem, see~\cite{DW}.
\item In the case of the symmetric random walk on the integers, one can prove this result with elementary tools, see~\cite{Desp}.  
\end{enumerate}
\end{rem} 

\subsection{Main result}
Recall that $\delta$ is a fixed minuscule weight and $P\left( \delta\right)$ is the finite set of weights of the minuscule representation $V\left( \delta\right) $. Let $\theta$ be in $\mathbb{R}_{>0}^d$ and $\left( p_s\right) _{s\in P\left( \delta\right) }$ the associated probability distribution on $P\left( \delta\right) $. Remind that $x$ is a solution in $\oplus_{i=1}^D\mathbb{R}_{>0}\varepsilon_i$ of System~\eqref{syst}.
 
From a probabilistic point of view, it is convenient to introduce a subset of $P^+$ which is more adapted to our situation. According to Proposition~\ref{combi}, the support of the conditioned random walk $\mathcal{W^+}$ defined in paragraph~\ref{llp} is the subset $P^{+,\delta}$ such that\[P^{+,\delta} =\left\lbrace \lambda\in P^+ :\exists n \geq 0 \enskip f_{\lambda/0,\delta}^n >0 \right\rbrace.\] Introduce also the sequence of subsets $\left( P_{n}^{+,\delta}\right) _{n\geq 0}$ defined by
\[P_{n}^{+,\delta} =\left\lbrace \lambda\in P^{+,\delta} : f_{\lambda/0,\delta}^n >0 \right\rbrace \qquad n\geq 0.\]
Observe that $P_{0}^{+,\delta}=\left\lbrace 0\right\rbrace$ and  \[P^{+,\delta}=\bigcup_{n\geq 0}P_{n}^{+,\delta}.\]
In general, we have $P^{+,\delta}\neq P_+$. With this notation, we have
   \[V\left(\delta \right)^{\otimes n} \simeq \bigoplus_{\lambda\in P^+} V\left( \lambda\right) ^{\oplus f_{\lambda/0,\delta}^n}=\bigoplus_{\lambda\in P_{n}^{+,\delta}} V\left( \lambda\right) ^{\oplus f_{\lambda/0,\delta}^n}\qquad n\geq 0.\]
By restricting the matrix $p$ to $P^{+,\delta}$ instead of $P^+$, we have a simple tool to derive our main result, see Lemma~\ref{compharm}.
   
The following identity can be proved by elementary computations on characters, see Proposition 5.3 in~\cite{LLP}.
\begin{equation}f_{\Lambda / \lambda,\delta}^{n}=\sum_{\kappa \in P_{n}^{+,\delta}} m_{\lambda,\kappa}^{\Lambda} f_{\kappa/0,\delta}^{n}\qquad \lambda,\Lambda \in P^{+},n\geq 0.\label{multi}\end{equation}

\begin{lem}
\label{recu}
We have\[\Lambda\in P_{n+1}^{+,\delta}\iff \exists\lambda\in P_{n}^{+,\delta} \enskip \lambda\leadsto\Lambda\qquad \Lambda\in P^+,\,n\geq 0. \]
\end{lem}
\textit{Proof.} On the one hand, we have\[V\left( \delta\right) ^{n+1}\simeq\bigoplus_{\Lambda\in P^+} V\left( \Lambda\right) ^{\oplus f_{\Lambda/0,\delta}^{n+1}}\]
and, in the other hand, we have\[V\left( \delta\right) ^{n+1}\simeq V\left( \delta\right) \otimes V\left( \delta\right) ^{\otimes n}\simeq V\left( \delta\right) \otimes \left(\bigoplus_{\lambda\in P^+} V\left( \lambda\right) ^{\oplus f_{\lambda/0,\delta}^{n}} \right).\]
Using Equality~\eqref{pieri}, this gives\[\begin{aligned}\sum_{\Lambda\in P^+}f_{\Lambda/0,\delta}^{n+1} s_{\Lambda}\left( x\right)&=s_{\delta}\left( x\right) \sum_{\lambda\in P^+}f_{\lambda/0,\delta}^n s_{\lambda}\left( x\right)=\sum_{\lambda\in P^+}f_{\lambda/0,\delta}^n s_{\lambda}\left( x\right)s_{\delta}\left( x\right) \\
&=\sum_{\lambda\in P^+}f_{\lambda/0,\delta}^n \left(\sum_{\lambda\leadsto\Lambda}s_{\Lambda}\left( x\right)  \right)=\sum_{\Lambda\in P^+}\left(\sum_{\lambda\leadsto\Lambda}f_{\lambda/0,\delta}^{n} \right)s_{\Lambda}\left( x\right). \end{aligned}\]
We deduce then the equality\[f_{\Lambda/0,\delta}^{n+1}=\sum_{\lambda\leadsto\Lambda} f_{\lambda/0,\delta}^n\]
which gives the expected result. $\Box$\\

Now we prove that the sequences\[\left(p_{n}^+\left( \lambda,\Lambda,i\right)  \right)_{n\geq 0}\qquad \lambda,\Lambda \in P^+,\,i\geq 0 \]
are convergent. Mimicking the previous section, Proposition~\ref{condifinite} leads us to introduce the two sequences of positive functions $\left(\psi_n=\psi_n\left( {\delta}\right)  \right)_{n\geq 1} $ and $\left(h_n=h_n \left( {\delta}\right) \right)_{n\geq 1}$ defined by
   \[\begin{array}{ccccc}
   \psi_{n} & : & P^{+,\delta} & \longrightarrow &  \mathbb{R}_{>0}\\
    & & \lambda & \longmapsto & \mathbb{P}_{\lambda}\left[ \mathcal{W}\in\overline{C}\right] =\mathbb{P}\left[\lambda+\mathcal{S}\left( \leq n\right) \in\overline{C}\right] 
   \end{array}\qquad
   \begin{array}{ccccc}
    h_{n} & : & P^{+,\delta} & \longrightarrow &  \mathbb{R}_{>0}\\
     & & \lambda & \longmapsto & \frac{\psi_n \left(\lambda \right) }{\psi_n \left( 0\right) }
    \end{array}\qquad n\geq 1.\] 

\begin{theo}
\label{harm}
Assume that the drift is zero (equivalently, $\theta=1^d$). The sequence $\left(h_n \right)_{n\geq 1} $ pointwise converges on $P^{+,\delta}$ to an harmonic function on $P^{+,\delta}$ for the restriction of the matrix $p$ to $P^{+,\delta}$.
\end{theo}
\textit{Proof.} Let $\lambda$ be in $P^{+,\delta}$ and $n$ a positive integer such that $n\geq2$. Remark that \[h_n \left(\lambda \right)=\frac{\psi_n \left(\lambda \right)}{\psi_n \left(0 \right)}=\frac{\mathbb{P}\left[\tau_\lambda>n \right] }{\mathbb{P}\left[\tau_0>n \right]}. \]
By Theorem~\ref{dw}, one has the pointwise convergence of the sequence $\left(h_n \right)_{n\geq 1}$ on $P^{+,\delta}$. Now we have \[\psi_n \left(\lambda \right)=\sum_{\lambda\leadsto\Lambda} p\left( {\lambda,\Lambda}\right) \psi_{n-1} \left(\Lambda \right)\qquad \lambda\in P^{+,\delta},\]
and so
\[h_n \left( \lambda\right) =\frac{\psi_{n-1} \left( 0\right) }{\psi_{n} \left( 0\right)}\sum_{\lambda\leadsto \Lambda} p{\left( \lambda,\Lambda\right) }h_{n-1} \left(\Lambda\right)\qquad \lambda\in P^{+,\delta}.\]
Theorem~\ref{dw} gives also \[\psi_{n-1} \left( 0\right)\underset{n\to \infty}{\sim}\psi_{n} \left( 0\right)\]
and the result follows. $\Box$\\

The following lemma helps us to identify the pointwise limit given by the previous theorem. This is an easy adaptation of a lemma called "Maximum principle" in~\cite{W}.
\begin{lem}
\label{compharm}
 Let $h_1,h_2$ be two harmonic functions on $P^{+,\delta}$ for the restriction of the matrix $p$ to $P^{+,\delta}$. If they coincide at $0$ and if $h_1\leq h_2$, then they coincide on $P^{+,\delta}$. 
 \end{lem}
 \textit{Proof.} For $h_1,h_2$ satisfying the assumptions of the proposition, consider the nonnegative function $h=h_2-h_1$. Both functions $h_1$ and $h_2$ being harmonic on $P^{+,\delta}$ for the restriction of $p$ to $P^{+,\delta}$, so is the function $h$. As $P_{0}^{+,\delta}=\left\lbrace 0\right\rbrace$, one has $h=0$ on $P_{0}^{+,\delta}$. Assume now that $h=0$ on $P_{n}^{+,\delta}$ for some nonnegative integer $n$ and suppose there exists $\Lambda_0$ in $P_{n+1}^{+,\delta}$ such that $h\left(\Lambda_0 \right)>0 $. Using Lemma~\ref{recu}, one can then find a weight $\lambda_0$ in $P_{n}^{+,\delta}$ such that $\lambda_0\leadsto\Lambda_0$. The function $h$ being harmonic and nonnegative on $P^{+,\delta}$, we obtain \[h\left(\lambda_0\right)=p\left(\lambda_0,\Lambda_0\right) h\left(\Lambda_0 \right)+ \sum_{\substack{\lambda_0\leadsto\Lambda\\ \Lambda\neq \Lambda_0}}  p\left( \lambda_0,\Lambda\right) h\left( \Lambda\right) >0.\] This is a contradiction. So $h$ is zero on $P_{n+1}^{+,\delta}$. By induction, $h$ is zero on $P^{+,\delta}$. $\Box$\\
 
 We already have an example of such positive harmonic function : by construction of the probability distribution $p$ and by Equality~\eqref{pieri}, the function
  \[\begin{array}{ccc} P^{+,\delta} & \longrightarrow & \mathbb{R}\\
  				\lambda & \longmapsto & x^{-\lambda} s_\lambda \left(x\right)
  				\end{array}\]  
 has the required property for every choice of $\theta$ in $\mathbb{R}_{>0}^d$. All that we have to do is to compare the pointwise limit given by Theorem~\ref{harm} to this candidate.
\begin{prop}
\label{comparaison}
We have \[h_n \left(\lambda \right) \leq x^{-\lambda} s_\lambda \left(x\right) \qquad \lambda\in P^{+,\delta},\,n\geq 1.\]
\end{prop} 
\textit{Proof.} When $n$ is a positive integer such that $n\geq 2$, one has
\[\begin{aligned}\psi_n \left(\lambda\right)&=\mathbb{P}\left[\bigcap_{i=1}^n \lambda+\mathcal{S}\left( i\right) \in \overline{C} \right]\\
&=\mathbb{P}\left[\bigcap_{i=1}^n  \lambda+\mathcal{S}\left( i\right) \in P^{+,\delta} \right]
&=\sum_{\Lambda\in P^{+,\delta} } \mathbb{P}\left[\left(\bigcap_{i=1}^{n-1} \lambda+\mathcal{S}\left( i\right) \in P^{+,\delta} \right) \cap\lambda+\mathcal{S}\left( n\right)  =\Lambda \right].\end{aligned}\]Given  $\Lambda$ in $P^{+,\delta}$, using~\eqref{centr}  and Proposition~\ref{combi}, one gets \[\begin{aligned}\mathbb{P}\left[\left(\bigcap_{i=1}^{n-1}\lambda+\mathcal{S}\left( i\right) \in P_{+}^\delta   \right) \cap \lambda+\mathcal{S}\left( n\right)=\Lambda \right]&=
\sum_{\substack{\pi\in \pi_{n}^+ \left( \lambda,\Lambda\right)\\\pi=\left( \pi_0,\pi_1,\ldots,\pi_n\right)}  }\mathbb{P}\left[\bigcap_{i=1}^n \lambda+\mathcal{S}\left( i\right) =\pi_i \right]\\&=\sum_{\pi\in \pi_{n}^+ \left( \lambda,\Lambda\right)} \frac{x^{\Lambda-\lambda}}{s_\delta \left( x\right) ^n}\\&=\mathrm{Card}\, \pi_{n}^+\left( \lambda,\Lambda\right) \frac{x^{\Lambda-\lambda}}{s_\delta \left( x\right) ^n}=f_{\Lambda/\lambda,\delta}^n \frac{x^{\Lambda-\lambda}}{s_\delta \left( x\right)^n }.\end{aligned}\]
This shows
\[\psi_n \left(\lambda\right)=\sum_{\Lambda\in P_{n}^{+,\delta}} f_{\Lambda/\lambda,\delta}^n \frac{x^{\Lambda-\lambda}}{s_\delta \left( x \right)^n}. \]
Now we use Identity~\eqref{multi} to write\[\psi_n \left(\lambda \right)=x^{-\lambda} \sum_{\Lambda\in P_{n}^{+,\delta}} \left(\sum_{\kappa\in P_{n}^{+,\delta}} f_{\kappa/0,\delta}^n m_{\lambda,\kappa}^\Lambda \right)  \frac{x^{\Lambda}}{s_\delta \left( x\right) ^n}=x^{-\lambda} \sum_{\Lambda\in P_{n}^{+,\delta}} \left(\sum_{\kappa\in P_{n}^{+,\delta}} f_{\kappa/0,\delta}^n m_{\kappa,\lambda}^\Lambda \right)  \frac{x^{\Lambda}}{s_\delta \left( x\right) ^n}.\]
Inequality~\eqref{cryst} gives then
\[\begin{aligned}\psi_n \left(\lambda \right)&\leq x^{-\lambda}\sum_{\Lambda\in P_{n}^{+,\delta}} \left(\sum_{\kappa\in P_{n}^{+,\delta}} f_{\kappa/0,\delta}^n K_{\lambda,\Lambda-\kappa} \right)  \frac{x ^{\Lambda}}{s_\delta \left( x\right) ^n}\\&=x^{-\lambda}\sum_{\kappa\in P_{n}^{+,\delta}} f_{\kappa/0,\delta}^n \left(\sum_{\Lambda\in P_{n}^{+,\delta}}  K_{\lambda,\Lambda-\kappa}x^{\Lambda-\kappa}\right) \frac{x^{\kappa}}{s_\delta \left( x\right) ^n}.\end{aligned}\]
Now remark that for any $\kappa$ in $P_{n}^{+,\delta}$, one has
\[\sum_{\Lambda\in P_{n}^{+,\delta}}  K_{\lambda,\Lambda-\kappa}x^{\Lambda-\kappa}\leq {\sum_{\gamma\in P} K_{\lambda,\gamma} x^{\gamma}=s_\lambda \left( x\right) }\]
and one gets
\[\psi_n \left(\lambda\right)\leq x^{-\lambda} s_\lambda \left( x \right)\sum_{\kappa\in P_{n}^{+,\delta}} f_{\kappa/0,\delta}^n \frac{x ^\kappa}{s_\delta \left( x \right) ^n}= x ^{-\lambda} s_\lambda \left( x \right) \psi_n \left( 0\right).\]$\Box$\\

We state now our main result. For the last time, recall that the drift is zero if and only if $\theta=1^d$, that is if and only if $p$ is the uniform distribution on $P\left( \delta\right)$ :\[p_s=\frac{1}{\dim  V\left( \delta\right) }\qquad s\in P\left( \delta\right).\]

\begin{theo}
\label{despax}
Assume that the drift of the random walk is zero. The sequence $\left(p_{n}^+ \left(\,.\,,\,.\,,0\right)\right)_{n\geq 0} $ pointwise converges on $P^{+,\delta} \times P^{+,\delta} $ and its limit $p^+=p^{+} \left(1^d\right) $ is given by\[p^+\left(\lambda,\Lambda\right)=p\left( 1^d\right)  \left( \lambda,\Lambda\right) \frac{s_\Lambda \left( 1^D\right) }{s_\lambda \left( 1^D\right) }=\begin{cases}\frac{1}{\dim V \left(\delta\right)  }\frac{\dim V  \left( \Lambda\right)   }{\dim V\left(\lambda\right) }&\text{if }\lambda\leadsto\Lambda\\
0&\text{otherwise}\end{cases}\qquad \lambda,\Lambda \in P^{+,\delta}. \]
By Theorem~\ref{llptheo}, we therefore have the pointwise convergence on $P^{+,\delta} \times P^{+,\delta} $ \[p^+ \left( \theta\right) \xrightarrow[\substack{\theta\to 1^d \\
\theta \in \left] 0,1\right[ ^d}]{} p^+\left(1^d\right) .\] 
\end{theo}
\textit{Proof.} In the zero drift case, Lemma~\ref{compharm} and Proposition~\ref{comparaison} show that the pointwise limit on $P^{+,\delta}$ of the sequence $\left( h_n\right) _{n\geq 1}$ is\[\lambda\mapsto s_\lambda \left(1^D\right)=\dim V\left( \lambda \right).\]Using Proposition~\ref{condifinite}, we obtain the expected result. $\Box$\\

To conclude, observe that the assumption on the drift is only required in Theorem~\ref{harm}. Thus, if Theorem~\ref{dw} is generalized in the near future to the case of a drift in the boundary of the cone, the generalization of our Theorem~\ref{despax} will automatically follow. 

\begin{conj}
	Assume the drift of the random walk is in $\partial{C}$. The sequence $\left(p_n^{+}\left(\,.\,,\,.\,,0\right)\right)_{n\geq 0} $ pointwise converges on $P^{+,\delta} \times P^{+,\delta} $ and its limit $p^+=p^+ \left( \theta\right) $ is given by\[p^+  \left(\lambda,\Lambda\right)=p  \left( \theta\right) \left( \lambda,\Lambda\right) \frac{x^{-\Lambda}s_\Lambda \left( x\right) }{x^{-\lambda}s_\lambda \left( x\right)}=\begin{cases}\frac{1}{s_\delta \left( x\right)   }\frac{s_\Lambda \left( x\right)    }{s_\lambda \left( x\right)  }&\text{if }\lambda\leadsto\Lambda\\
	0&\text{otherwise}\end{cases}\qquad \lambda,\Lambda \in P^{+,\delta}. \]
\end{conj}

\textbf{Acknowlegments.} The author would like to deeply acknowledge O. Durieu, E. Journé, C. Lecouvey, E. Lesigne, M. Peigné, K. Raschel for helpfull discussions and support during the writing of this article.

\end{document}